\documentclass{article}
\usepackage{amsmath}
\usepackage{amssymb}
\usepackage{verbatim}
\usepackage{enumerate}

\newcommand{\pmax}{\mathbb{P}_{max}}

\newcommand{\psmax}{\mathbb{P}^{*}_{max}}

\newcommand{\calp}{\mathcal{P}}

\newtheorem{df}{Definition}[section]
\newtheorem{thrm}[df]{Theorem}
\newtheorem{lem}[df]{Lemma}

\begin{document}
\pagenumbering{arabic}
\title{Showing OCA in $\pmax$-style extensions}
\author{Paul Larson\thanks{This material is based in part upon work 
supported by the North Atlantic Treaty Organization under a grant awarded 
in 1998.}}
\maketitle

\begin{abstract} We give a proof that OCA holds in the $\pmax$ extension 
of $L(\mathbb{R})$. The proof is general enough to be adapted 
to most $\pmax$ variations.
\end{abstract}

\section{Introduction}

This paper is an essentially expository presentation of a proof that OCA
holds in the $\pmax$ extension of $L(\mathbb{R})$, assuming
AD$^{L(\mathbb{R})}$. The proof consists of putting together standard
facts about $\pmax$ and OCA; essentially none of the ideas here are 
due to the author. One purpose in writing this is to make public the 
details of a remark made in $\cite{LT}$. Another is to show that OCA
is achieved fairly easily in $\pmax$-style extensions. One possible
application of this fact is that if it turns out to be possible to have
$\pmax$ variations with $c > \omega_{2}$, OCA should hold in such 
an extension.

Instead of $\pmax$, however, we will work with a variation called $\psmax$. 
The 
main difference between the two is that $\pmax$ conditions are individual 
models, while $\psmax$ conditions are sequences of models. $\psmax$ is 
really just the limit points of $\pmax$. It is a standand fact $\cite{W}$ 
(under AD$^{L(\mathbb{R})}$) that every $\pmax$ extension is $\psmax$ 
extension and vice versa.

\section{Iterable structures}\label{itst}

The following is the definition of iterability for sequences of models.

\begin{df}\label{itdef}([3]) Suppose $\langle N_{k} : k < \omega\rangle$
is a
countable
sequence
such that for each $k,$ $N_{k}$ is a countable transitive model of
ZFC$^{*}$ and
such that for all $k$, $N_{k} \in N_{k+1}$ and $\omega_{1}^{N_{k}} =
\omega_{1}^{N_{k+1}}$. An iteration of $\langle N_{k} : k < \omega\rangle$
is a
sequence

$$\langle\langle N^{\beta}_{k} : k < \omega \rangle , G_{\alpha} ,
j_{\alpha , \beta} : \alpha < \beta < \gamma \rangle$$
such that for all $\alpha < \beta < \gamma$ the following hold.

\begin{enumerate}

\item $j_{\alpha , \beta} : \cup \{ N^{\alpha}_{k} \mid k < \omega \}
\rightarrow \cup \{ N^{\beta}_{k} \mid k < \omega \}$ is a commuting
family
of $\Sigma_{0}$ elementary embeddings.

\item For all $k < \omega $, $G_{\beta} \cap N^{\beta}_{k}$ is an 
$N^{\beta}_{k}$-normal ultrafilter on
$(\calp(\omega_{1}))^{N^{\beta}_{k}}$.

\item If $\beta + 1 < \gamma$ then $N^{\beta + 1}_{k}$ is the $\cup \{
N^{\beta}_{k} \mid k < \omega \}$-ultrapower of $N^{\beta}_{k}$ by
$G_{\beta}$ and $j_{\beta , \beta + 1} : \cup \{ N^{\beta}_{k} \mid k <
\omega \}
\rightarrow \cup \{ N^{\beta + 1}_{k} \mid k < \omega \}$ is the induced
$\Sigma_{0}$ elementary embedding.

\item For each $\beta < \gamma$ if $\beta$ is a limit ordinal then for
every $k < \omega$, $N^{\beta}_{k}$ is the direct limit of $\{
N^{\alpha}_{k} \mid \alpha < \beta \}$ and for all $\alpha < \beta$,
$j_{\alpha , \beta }$ is the induced
$\Sigma_{0}$ elementary embedding.

\end{enumerate}

If $\gamma$ is a limit ordinal then $\gamma$ is the length of the 
iteration, otherwise the length of the iteration is $\delta$ where
$\delta + 1 = \gamma$.

A sequence $\langle N_{k}^{*} : k < \omega\rangle$ is an iterate of   
$\langle
N_{k} : k < \omega\rangle$ if it occurs in an iteration of $\langle N_{k}
: k <
\omega\rangle$.

The sequence $\langle N_{k} : k < \omega\rangle$ is iterable if every
iterate
of it is
well founded.

If $B \subset \mathbb{R}$, then  $\langle N_{k} : k < \omega\rangle$ is 
$B$-iterable if it is iterable, and if 
for every iterate $\langle N_{k}^{*} : k < \omega\rangle$ of $\langle N_{k} 
: k < \omega\rangle$, $j(B \cap N_{0}) = B \cap N_{0}^{*}$, where $j$ is the 
induced embedding and $j(B \cap N_{0})$ is defined to be $\cup \{j(a) : a 
\in N_{0} \text{ and } a \subset B \}$.

\end{df}

The following lemma is the main tool for verifying that the sequences as 
above are iterable.

\begin{lem}\label{seqit}([3]) Suppose
$$\langle N_{k} : k < \omega \rangle$$
is a countable sequence such that for each $k$, $N_{k}$ is a countable
transitive model of ZFC$^{*}$ and such that for all $k$, 
$$N_{k} \in N_{k+1}$$
and
$$(\omega_{1})^{N_{k}} = (\omega_{1})^{N_{k+1}}.$$
Suppose that for all $k < \omega$
\begin{enumerate}[(i)]
\item if $C \in N_{k}$ is closed and unbounded in $\omega_{1}^{N_{0}}$,
then there exists $D \in N_{k+1}$ such that $D \subset C$, $D$ is closed
and unbounded in $C$, and
$$D \in L[x]$$
for some $x \in \mathbb{R} \cap N_{k+1}$.
\item for all $x \in \mathbb{R} \cap N_{k}$, $x^{\#} \in N_{k+1}$.
\item for all $k < \omega$,
$$|N_{k}|^{N_{k+1}} = (\omega_{1})^{N_{0}}.$$
\end{enumerate}
Then the sequence $\langle N_{k} : k < \omega\rangle$ is iterable.
\end{lem}

We quote a lemma from [3] showing that under certain 
circumstances the ultrafilter needed to iterate a given sequence exists.

\begin{lem} \label{agen} ([3]) Suppose that
$$\langle N_{k} : k < \omega \rangle$$
is a sequence of countable transitive sets such that for all $k < \omega$,
$N_{k} \in N_{k+1}$, 
$$N_{k} \models \text{ZFC}^{*},$$
and
$$N_{k} \cap (I_{NS})^{N_{k+1}} = N_{k+1} \cap (I_{NS})^{N_{k+2}}.$$
Suppose that $k \in \omega$ and that
$$a \in (\calp(\omega_{1}))^{N_{k}}\setminus (I_{NS})^{N_{k+1}}.$$
Then there exists
$$G \subset \cup \{(\calp(\omega_{1}))^{N_{i}} \mid i < \omega \}$$
such that $a \in G$ and such that for all $i < \omega$, $G \cap N_{i}$ is
a uniform $N_{i}$-normal ultrafilter.
\end{lem}




The sequences of models in $\mathbb{P}_{max}$ variations satisfy a 
variation of $\psi_{AC}$.

\begin{df}\label{pac}([3]) $\psi^{*}_{AC}$: Suppose that $\langle
S_{\alpha} : \alpha < \omega_{1} \rangle$ and $\langle T_{\alpha} : \alpha
< \omega_{1} \rangle$ are each sequences of stationary, costationary sets.
Then there exists a sequence $\langle \delta_{\alpha} : \alpha <
\omega_{1} \rangle$ of ordinals less than $\omega_{2}$ such that for each
$\alpha < \omega_{1}$ there exists a bijection
$$\pi : \omega_{1} \rightarrow \delta_{\alpha},$$
and a closed unbounded set $C \subset \omega_{1}$ such that
$$\{\eta < \omega_{1} \mid \text{ }o.t.(\pi [\eta]) \in T_{\alpha} \} \cap
C 
= S_{\alpha} \cap C.$$
\end{df}

The reason for this variation is that our conditions are sequences of
models, and 
iterates of sequences modeling $\psi^{*}_{AC}$ model $\psi^{*}_{AC}$. This
isn't so
for $\psi_{AC}$.

\section{$\psmax$}

\begin{df}([3]) $\psmax$ is the set of pairs $( \langle
M_{k} : k
< \omega \rangle, a)$ such that the following hold.

\begin{enumerate}

\item $a \in M_{0}, a \subset \omega_{1}^{M_{0}}$, and
$\omega_{1}^{M_{0}}=\omega_{1}^{L[a,x]}$ for some $x \in \mathbb{R}
\cap M_0$.

\item Each $M_{k}$ is a countable transitive model of ZFC$^{*}$.

\item $M_{k}\in M_{k+1}, \omega_{1}^{M_{k}}=\omega_{1}^{M_{k+1}}$.

\item $(I_{NS})^{M_{k+1}}\cap M_{k}=(I_{NS})^{M_{k+2}}\cap M_{k}$.  

\item $\cup\{M_{k}:k<\omega\}\models\psi^{*}_{AC}$.

\item $\langle M_{k} \mid k <\omega\rangle$ is iterable.

\item $\exists X \in M_0$ such that
$X \subset \calp(\omega_{1})^{M_{0}}\setminus I^{M_{1}}_{NS}$, such
that
$M_{0}\models ``|X| = \omega_{1},"$ and such that for all $A,B
\in X$, if $A\neq B$ then $A \cap B \in I ^{M_{0}}_{NS}$.
\end{enumerate}

The ordering on $\psmax$ is as follows.

$$(\langle N_{k} :k <\omega\rangle ,  b) <
(\langle M_{k} :k <\omega\rangle, a)$$ if
$\langle M_{k} :k <\omega\rangle \in N_{0},
\langle M_{k} :k <\omega\rangle$ is hereditarily countable in $N_{0}$
and there exists an iteration $$j: \langle M_{k} : k <\omega\rangle
\rightarrow \langle M_{k}^{*} :k <\omega\rangle$$ such that:

\begin{enumerate}
\item $j(a) = b$,

\item $\langle M_{k}^{*} :k <\omega\rangle \in N_{0}$ and $j \in N_{0}$,

\item $j(I^{M_{k+1}}_{NS})\cap M^{*}_{k} =(I_{NS})^{N_{1}} \cap M^{*}_{k}$ 
for all $k <\omega$,

\end{enumerate}

\end{df}

The following lemma follows from the fact that $\psmax$ conditions model
$\psi^{*}_{AC}$. The analogous lemma holds in the other
$\mathbb{P}_{max}$ variations whose conditions are sequences of models.

\begin{lem}\label{auniq}([3]) Suppose that $(\langle M_{k} \mid k <
\omega \rangle, a) \in \psmax$. Suppose that
$$j_{1} : \langle M_{k} \mid k < \omega \rangle\rightarrow \langle
M^{1}_{k} \mid k < \omega \rangle$$
and
$$j_{2} : \langle M_{k} \mid k < \omega \rangle\rightarrow \langle
M^{2}_{k} \mid k < \omega \rangle$$
are well founded iterations such that $j_{1}(a) = j_{2}(a)$.

Then
$$\langle M^{1}_{k} \mid k < \omega \rangle=\langle M^{2}_{k} \mid k <
\omega \rangle$$
and $j_{1} = j_{2}$.
\end{lem}                 

Since the order on $\psmax$ is determined by the
existence of elementary embeddings, each condition $(\langle M_{k} : k <
\omega \rangle, a)$ in the
generic is iterated $\omega_1$ times through the conditions below it in
the generic.  In fact, by Lemma $\ref{auniq}$, each $(\langle M_{k} : k <
\omega \rangle, a)$ in the generic is uniquely iterated
into the extension to
a structure
$\langle\langle M^*_{k} : k <\omega\rangle, a_{G}\rangle$, where
$$a_{G}=\cup\{ a \mid \exists (\langle M_{k} : k < \omega \rangle, a) \in
G \},$$ for generic $G$. The following definitions apply to all $\mathbb{P}_{max}$ variations.

\begin{df}($\cite{W}$) A filter $G\subset\psmax$ is semi-generic if for all 
$\alpha <
\omega_{1}$ there exists a condition $\langle M_{k} : k < \omega_{1}  
\rangle \in G$ such that $\alpha < \omega_{1}^{M_{0}}.$

$A_{G} = \cup\{ a \mid \exists (\langle M_{k} : k < \omega \rangle, a)
\in G \}$. 
$$\calp(\omega_1)_G = \bigcup \{\calp(\omega_1)^{M^*_{0}}
\mid (\langle M_{k} : k < \omega \rangle, a) \in G\},$$ and
$$I_{G} = \cup \{I_{NS}^{M^{*}_{1}} \cap M^{*}_{0} \mid   
(\langle M_{k} : k < \omega \rangle, a) \in G \},$$ where for $(\langle 
M_{k} : k < \omega \rangle, a) \in G$, $\langle M^{*}_{k} : k < 
\omega\rangle$ is the iterate of $\langle M_{k} : k < \omega \rangle$ by 
the unique iteration of $\langle M_{k} : k < \omega \rangle$ that sends 
$a$ to $A_{G}$. 

\end{df}

The following lemma summarizes the basic analysis of $\psmax$.

\begin{thrm}([3]) Assume $\text{AD}^{L(\mathbb{R})}$.
Then $\psmax$ is $\omega$-closed and homogeneous.

Suppose G
$\subset \psmax$ is $L(\mathbb{R})$-generic.
Then
$$L(\mathbb{R})[G] \models \omega_{1}-\text{DC}$$
and in $L(\mathbb{R})[G]$:
\begin{enumerate}
\item $\calp(\omega_1)_G = \calp(\omega_1)$.
\item $I_{G}$ is a normal saturated ideal;
\item $I_{G}$ is the nonstationary ideal.
\end{enumerate}
\end{thrm}

The following theorem is implicit in [3], and is the key to the
main theorem.

\begin{thrm} \label{condex}([3]) Assume AD holds in $L(\mathbb{R})$.
Then for each set
$A\subset\mathbb{R}$ with $A\in L(\mathbb{R})$ there is a condition
$(\langle M_{k} : k < \omega \rangle, a) \in \psmax$ such
that each $M_{k}$ is a model of ZFC, and such that 
for all stationary set preserving forcings $P$ in $M_{0}$, if
$G \subset P$ is $M_{k}$-generic for all $k < \omega$, then 
$(\langle M_{k}[G] : k < \omega \rangle, a) \in \psmax$ and
\begin{enumerate}
 
\item $A \cap M_{0}[G] \in M_{0}[G],$

\item $\langle H(\omega_{1})^{M_{0}[G]}, A \cap M_{0}[G] \rangle 
\prec \langle H(\omega_{1}), A
\rangle,$   

\item $\langle M_{k}[G] : k < \omega \rangle$ is $A$-iterable.
\end{enumerate}
Moreover, the set of such conditions is dense in
$\psmax$.
\end{thrm}

\section{OCA}

\begin{df}The Open Coloring Axiom (OCA) is the statement that if $O
\subset \mathbb{R} \times \mathbb{R}$ is open and symmetric, and $A
\subset \mathbb{R}$, then either there is an uncountable set $B \subset
A$ such that $[B]^{2} \subset O$, or $A$ is the union of countably many  
sets $\langle C_{n} \mid n \in \omega \rangle$ such that each
$[C_{n}]^{2} \cap O = \emptyset$.
\end{df}

The relevant theorem about OCA is the following.

\begin{thrm}\label{foca}($\cite{T}$) Say that 
$A \subset \mathbb{R} \times \mathbb{R}$
is open and symmetric, and that $X \subset \mathbb{R}$
is not the union of $\omega$ many sets $Y$ such that
$[Y]^{2} \cap A = \emptyset$. Then there is an
$X' \subset X$ of cardinality $2^{\omega}$ such that
the partial order consisting of finite $a \subset X'$ 
such that $[a]^{2} \subset A$, ordered by inclusion,
is $<2^{\omega}$-c.c.
\end{thrm}

The importance of this theorem is that if CH holds then 
there is a c.c.c. forcing to obtain a set witnessing any
given instance of OCA. 

\begin{thrm} Assume $AD^{L(\mathbb{R})}$, and let $G \subset \psmax$ be 
$L(\mathbb{R})$-generic. Then OCA holds in $L(\mathbb{R})[G]$.
\end{thrm}

Proof: Let $x$ be a real coding an open, symmetric subset $O$ of 
$\mathbb{R} \times 
\mathbb{R}$ and let $\tau \in L(\mathbb{R})$ be a $\psmax$-name for a set 
of reals. Let $( \langle N_{k} : k < \omega \rangle, a) \in \psmax$ force 
that the realization of $\tau$ is not 
the union of $\omega$ many sets $A_{i}$ such that for each $i$ and each 
$x, y \in A_{i}$, $(x, y) \not\in O$. Let $z$ be a real coding $( \langle 
N_{k} : k < \omega \rangle, a)$, let $B$ be a 
set of reals 
coding $\tau$ and let $A = B \times \{z\}$. For this set $A$, let
$(\langle M_{k} : k < \omega \rangle, b )$ be a $\psmax$ 
condition as given by Lemma 
$\ref{condex}$. We may assume by forcing over 
$\bigcup\{M_{k} : k < \omega \}$ if necessary that CH holds
in $M_{0}$.

In $M_{0}$, build a decreasing $\omega_{1}$-sequence $\langle p_{\alpha} 
= (\langle M^{\alpha}_{k} : k < \omega \rangle, a^{\alpha}) 
: \alpha < \omega_{1}^{M_{0}} \rangle$ of $\pmax$ conditions as follows. 
Let $\langle S^{\alpha}_{i} : i < \omega,\text{ }\alpha < 
\omega_{1}^{M_{0}} \rangle \in M_{0}$ be a set of mutually disjoint 
subsets of $\omega_{1}^{M_{0}}$ which are stationary in $M_{1}$, and let 
$\langle A^{\alpha}_{i} : i < \omega,\text{ }\alpha < \omega_{1}^{M_{0}} 
\rangle \in M_{0}$ be a listing of all the countable sequences of closed 
sets of reals in $M_{0}$. As we build our sequence, let $j_{\alpha, 
\beta}$ be the embedding witnessing that $p_{\beta} < p_{\alpha}$. 
Further, as we build our sequence, build $\langle 
B^{\alpha}_{i} : i < \omega,\text{ } \alpha < \omega_{1}^{M_{0}} \rangle$ 
to enumerate $ \bigcup\{ \mathcal{P}(\omega_{1})^{M^{\alpha}_{k}} 
\setminus (I_{NS}^{M^{\alpha}_{k + 1}}) : k < \omega, \alpha < 
\omega^{M_{0}}_{1} \}$.

Let $( \langle N_{k}: k < \omega \rangle, a) = p_{0}$. Given 
$p_{\alpha}$, the $\alpha$-th member of our sequence, we choose 
$p_{\alpha + 1}$ as follows. Since $H(\omega_{1})^{M_{0}}$ is 
elementary in $H(\omega_{1})$ with the predicate $A$, for each condition 
$q$ in $(\psmax)^{M_{0}}$ $(= \psmax \cap M_{0})$, and each countable 
sequence $\langle A_{i} : i < \omega \rangle$ of 
closed subsets of $\mathbb{R} \times \mathbb{R}$ in 
$M_{0}$, there is a condition $q' < q$ such that either $q'$ forces some 
real 
to be in $\tau$ which is not in the union of the closed sets, or $q'$ 
forces two reals $x$ and $y$ such that $(x, y) \in O$ to be in the same 
$A_{i}$. Given $p_{\alpha}$ as $q$ and $\langle A^{\alpha}_{i} : i < 
\omega \rangle$ as $\langle A_{i} : i < \omega \rangle$, let $p_{\alpha + 
1}$ be the corresponding $q'$.

The argument at limit stages is standard. To pick $p_{\alpha}$, we apply
the proof of the $\omega$-closure of $\psmax$. Let $\langle \alpha_{i} : i
< \omega \rangle$ be an increasing cofinal sequence below $\alpha$. Let
$j_{\alpha_{i} \infty}$ be the composition of all the $j_{\alpha{e}
\alpha_{e+1}}$ for $e \geq i$. Then $\langle M^{*}_{k} : k < \omega 
\rangle = \langle j_{\alpha_{k}
\infty}(M^{\alpha_{k}}_{k}) : k < \omega \rangle$ is an iterable 
sequence, by Lemma $\ref{seqit}$. Applying Theorem $\ref{condex}$ in 
$M_{0}$, pick a $\psmax$ condition $( \langle
\bar{M}_{k} : k < \omega \rangle, d)$ containing a real coding $\langle 
M^{*}_{k} : k < \omega \rangle$.
Then in $\bar{M}_{0}$, we can build an iteration
$j$ of $\langle M^{*}_{k} : k < \omega \rangle$ of length
$\omega_{1}^{\bar{M}_{0}}$ such that $\omega_{1}^{M^{*}_{0}}
\in j(j_{\alpha_{i}, \infty}(j_{\gamma, \alpha_{i}}(B^{\gamma}_{k})))$ if
$\omega_{1}^{M^{*}_{0}} \in S^{\gamma}_{k}$, for $k <
\omega$ and $\gamma < \omega_{1}^{M_{0}}$, and such that for each $k < 
\omega$, $j(\mathcal{P}(\omega_{1}^{M{*}_{k}}) \cap I_{NS}^{M^{*}_{k+1}}) 
= j(\mathcal{P}(\omega_{1}^{M{*}_{k}})) \cap I_{NS}^{\bar{M}_{1}}$. Then 
letting $a^{\alpha} = j(j_{\alpha_{0} \infty}(a^{\alpha_{0}}))$ and 
$j_{\alpha_{i} \alpha} = j \circ j_{\alpha_{i} \infty}$, we let 
$p_{\alpha} = ( \langle \bar{M}_{k} : k < \omega \rangle, a^{\alpha})$.

By Theorem $\ref{foca}$, there is a c.c.c. forcing in $M_{0}$ to get a 
$O$-homogeneous set of cardinality $\omega_{1}^{M_{0}}$
contained in the set of reals $y$  for which some $p_{\alpha}$ forces 
$y$ to be in $\tau$. Let $X$ be 
$\bigcup\{M_{k} : k < \omega \}$-generic for this forcing, and let $a$ be 
the union of the $a^{\alpha}$. Then $( \langle M_{k}[X] : k < \omega 
\rangle, a)$ is a $\psmax$ condition below every $p_{\alpha}$, and
by Theorem $\ref{condex}$, $( \langle M_{k}[X] : k < \omega 
\rangle, a)$ is $A$-iterable.
Then for every $\psmax$ generic $G$ with 
$(\langle M_{k}[X] : k < \omega \rangle, a) \in G$, if $j$ is the 
unique iteration sending $\langle M_{k}[X] : k < \omega \rangle$ 
through the generic, $j(X) \subset \tau_{G}$ is a witness for the
given instance of OCA. 
$\blacksquare$ 
\vspace{\baselineskip}

\noindent{Department of Computer\\
and Systems Engineering\\ Kobe University\\ Kobe 657-8501\\ 
Japan\\}
 
\noindent{larson\@@alan.scitec.kobe-u.ac.jp}

\end{document}